\documentclass{amsart}
\input cyracc.def

\usepackage{amssymb}
\usepackage{amsthm}

\newcommand{\sA}{\mathcal{A}}
\newcommand{\sM}{\mathcal{M}}
\newcommand{\sP}{\mathcal{P}}
\newcommand{\fN}{\mathfrak{N}}
\newcommand{\sN}{\mathcal{N}}
\newcommand{\sB}{\mathcal{B}}
\newcommand{\sT}{\mathcal{T}}

\newcommand{\real}{\mathbb{R}}

\newcommand{\Z}{\mathbb{Z}}

\newtheorem{thm}{Theorem}[section]
\newtheorem{lem}[thm]{Lemma}
\newtheorem{cor}[thm]{Corollary}

\newtheorem{prop}[thm]{Proposition}
\newtheorem{dfn}[thm]{Definition}

\title{On three notions of effective computation over $\real$}
\author{Wesley Calvert}
\address{Department of Mathematics \& Statistics\\ Faculty Hall 6C\\ Murray State
  University\\ Murray, Kentucky 42071}
\email{wesley.calvert@murraystate.edu}
\date{\today}
\thanks{The author is a member of the Department of Mathematics \&
  Statistics at Murray State University.  Correspondence and proofs
  should be sent to the author at {\tt
    wesley.calvert@murraystate.edu}.  The author may be reached by
  telephone at 1-270-809-2503, but callers should be aware that this
  location is in US Central time (GMT-6:00).}

\begin{document}


\begin{abstract}
We compare three notions of effectiveness on uncountable structures.
The first notion is that of a $\real$-computable structure, based on a
model of computation proposed by Blum, Shub, and Smale, which uses
full-precision real arithmetic.  The second notion is that of an
$F$-parameterizable structure, defined by Morozov and based on
Mal'tsev's notion of a constructive structure.  The third is
$\Sigma$-definability over $HF(\real)$, defined by Ershov as a
generalization of the observation that the computably enumerable sets
are exactly those $\Sigma_1$-definable in $HF(\mathbb{N})$.

We show that every $\real$-computable structure has an
$F$-parameterization, but that the expansion of the real field by the
exponential function is $F$-parameterizable but not
$\real$-computable.  We also show that the structures with $\real$-computable
copies are exactly the structures with copies $\Sigma$-definable over
$HF(\real)$.  One consequence of this equivalence is a method of
approximating certain $\real$-computable structures by Turing
computable structures.\end{abstract}

\maketitle

\newpage

\section{Introduction}

In the mid-twentieth century, Russian and Western mathematicians gave
definitions clarifying which algebraic structures could be considered
computable.  The Western approach focused on the effectiveness of the
atomic diagram of the structure:

\begin{dfn}[\cite{frsh, rabin}]\label{tcompstr} Let $\sA$ be a
  structure in a finite language.  We say that $\sA$ is
  \emph{computable} if and only if the atomic diagram of $\sA$ is
  Turing computable. \end{dfn}

Of course, many equivalent definitions could be given for Turing
computation, and the definition is often broadened to allow for a
\emph{computable} language, rather than a merely finite one.  Also,
this notion is not isomorphism invariant.  Many structures are not
themselves computable, but they have isomorphic copies that are.  The
key Russian definition was equivalent, but had a different focus:

\begin{dfn}[\cite{maltsev}]\label{constr} Let $\sA$ be a structure in a finite
  language.  Let $\nu :D \to \sA$ be surjective, for some $D \subseteq
  \omega$.  Then the \emph{enumerated structure} $(\sA, \nu)$ is said to be
  \emph{constructive} if and only if the following hold:
\begin{enumerate}
\item The set $\nu^{-1}$ is Turing computable.
\item For each predicate $P$ in the language (including the graphs of
  any functions in the language), the set $\nu^{-1}(P)$ is
  Turing computable.
\end{enumerate}
\end{dfn}

A key restriction in each definition is that the structure $\sA$ must be
countable in order to satisfy either.  In the last few years,
attention from various directions has turned to the question of
effectiveness for uncountable structures.  In the most technical
sense, the additive group of real numbers fails to be effective
(computable).  However, in any sense but the most technical, addition
of real numbers is usually thought to be quite effective.

There have been many attempts to formalize this intuition.  One model,
proposed by the present author in joint work with Porter \cite{cp1}, is
to generalize Definition \ref{tcompstr} by replacing Turing
computability with a model of computation proposed by Blum, Shub, and
Smale, in which full-precision arithmetic is axiomatically
``computable'' (a more careful definition will be given
later).  A structure satisfying this hypothesis is said to be
$\real$-computable.  Another, proposed by Morozov \cite{fparam}, generalizes
Definition \ref{constr} by replacing $\nu$ with a map from $\sA$ into
Baire space, where ``Turing computable'' is replaced by ``analytic''
(again, a formal definition will follow).  Such a structure is said to
be $F$-parameterizable.  A third, proposed by Ershov
\cite{ershovopred}, generalizes a classical theorem that computably enumerable sets
are exactly those $\Sigma_1$-definable in $HF(\mathbb{N})$ by
replacing $\mathbb{N}$ with the real field.  This kind of structure is
said to be $\Sigma$-definable in $HF(\real)$.  The purpose of the present paper is to
compare these three models for effectiveness on uncountable structures.

Other formalizations of effectiveness for uncountable structures
include: the notion of ``local computability,'' studied by R.\ Miller
and others \cite{rmiller, rmiller2}; infinite time machines, studied by Hamkins, R.\
Miller, Seabold, and Warner \cite{itcmt}; B\"uchi automata
\cite{montalbanbuchi}; and Borel structures, studied by H.\ Friedman and
C.\ Steinhorn \cite{steinhornborel1}.

In Section \ref{dfns} we will give the definitions of the three
classes of structures.  In Section
\ref{compare}, we will show that the class of $F$-parameterizable
structures properly contains that of $\real$-computable functions.  In
Section \ref{rcmtsig}, we will show that the classes of
$\real$-computable and $\Sigma$-definable in $HF(\real)$ structures
are equivalent.

\section{Definitions}\label{dfns}

\subsection{$\real$-Computability}
The definition of $\real$-computable structures comes from
\cite{cp1}, and the BSS model of computation is explained in detail
there.  The idea of this model is that full precision real number
arithmetic is axiomatically effective.  We give here an outline of the
relevant definitions.

The definition of a BSS machine comes from \cite{bss}, and the
concept is more fully described in \cite{bcss}.  Let $\real^\infty$ be the set of finite sequences of
elements from $\real$, and $\real_{\infty}$ the bi-infinite direct sum
\[\bigoplus\limits_{i \in \Z} \real.\]

\begin{dfn} A machine $M$ over $\real$ is a finite connected directed
  graph, containing five types of nodes: input, computation, branch,
  shift, and output, with the following properties:
\begin{enumerate}
\item The unique input node has no incoming edges and only one
  outgoing edge.
\item Each computation and shift node has exactly one output edge and
  possibly several input branches.
\item Each output node has no output edges and possibly several input
  edges.
\item Each branch node $\eta$ has exactly two output edges (labeled $0_\eta$ and $1_\eta$) and possibly
  several input edges.
\item Associated with the input node is a linear map $g_I:\real^\infty \to
  \real_\infty$.
\item Associated with each computation node $\eta$ is a rational
  function $g_\eta: \real_\infty \to \real_\infty$.
\item Associated with each branch node $\eta$ is a polynomial function
  $h_\eta : \real_\infty \to \real$.
\item Associated with each shift node is a map $\sigma_\eta \in
  \{\sigma_l, \sigma_r\}$, where $\sigma_l(x)_i = x_{i+1}$
  and $\sigma_r(x)_i = x_{i-1}$.
\item Associated with each output node $\eta$ is a linear map $O_\eta:
  \real_\infty \to \real^{\infty}$.
\end{enumerate}
\end{dfn}

A machine may be understood to compute a function in the following
way:
\begin{dfn} Let $M$ be a machine over $\real$.
\begin{enumerate}
\item A \emph{path} through $M$ is a sequence of nodes $(\eta_i)_{i
  =0}^n$ where $\eta_0$ is the input node, $\eta_n$ is an output node,
  and for each $i$, we have an edge from $\eta_i$ to $\eta_{i+1}$.
\item A \emph{computation} on $M$ is a sequence of pairs
  $\left((\eta_i, x_i)\right)_{i = 0}^n$ with a number $x_{n+1}$, where $(\eta_i)_{i
  =0}^n$ is a path through $M$, where $x_0 \in R^\infty$, and where,
  for each $i$, the following hold:
\begin{enumerate}
\item If $\eta_i$ is an input node, $x_{i+1} = g_I(x_i)$.
\item If $\eta_i$ is a computation node, $x_{i+1} = g_{\eta_i}(x_i)$.
\item If $\eta_i$ is a branch node, $x_{i+1} = x_i$ and $\eta_{i+1}$
  determined by $h_{\eta_i}$ so that if $h_{\eta_i}(x_i) \geq 0$, then
  $\eta_{i+1}$ is connected to $\eta_i$ by $1_{\eta_i}$ and if $h_{\eta_i}(x_i) < 0$, then
  $\eta_{i+1}$ is connected to $\eta_i$ by $0_{\eta_i}$.  (Note that
  in all other cases, $\eta_{i+1}$ is uniquely determined by the
  definition of path.)
\item If $\eta_i$ is a shift node, $x_{i+1} = \sigma_{\eta_i}(x_i)$
\item If $\eta_i$ is an output node, $x_{i+1} = O_{\eta_i}(x_i)$.
\end{enumerate}
\end{enumerate}
\end{dfn}

The proof of the following lemma is an obvious from the definitions.

\begin{lem} Given a machine $M$ and an element $z \in \real^{\infty}$,
  there is at most one computation on $M$ with $x_0 = z$. \end{lem}

\begin{dfn} The function $\varphi_M : \real^\infty \to R^\infty$ is
  defined in the following way:  For each $z \in \real^\infty$, let
  $\varphi_M(z)$ be $x_{n+1}$, where $\left(\left((\eta_i,
  x_i)\right)_{i = 0}^n, x_{n+1}\right)$ is the unique computation, if
  any, where $x_0 = z$.  If there is no such computation, then
  $\varphi_M$ is undefined on $z$.\end{dfn}

Since a machine is a finite object, involving finitely many real
numbers as parameters, it may be coded by a member of $\real^\infty$.

\begin{dfn}If $\sigma$ is a code for $M$, we define $\varphi_\sigma =
  \varphi_M$.\end{dfn}

We can now say that a set is $\real$-computable if and only if its
characteristic function is $\varphi_M$ for some $M$.  Now by
identifying a structure with its atomic diagram in the usual way, we
define a set to be $\real$-computable if and only if its atomic
diagram is $\real$-computable.

\subsection{$F$-Parameterizability}
In \cite{fparam}, Morozov introduced a concept that he called
\emph{$F$-parameterizability} in order to understand the elementary
substructure relation on both automorphism groups and the structure of
hereditarily finite sets over a given structure.  In a talk at
Stanford University, though, he identified this notion as one ``which
generalizes the notion of computable'' \cite{morozovstanford}.

\begin{dfn}[\cite{fparam}]\label{df:fparam} Let $\sM$ be a structure
  in a finite relational language $\left(P_n^{k_n})_{n \leq k}\right)$.  We
  say that $\sM$ is $F$-parameterizable if and only if there is an
  injection $\xi : \sM \to \omega^\omega$ with the following
  properties:
\begin{enumerate}
\item The image of $\xi$ is analytic in the Baire space, and
\item For each $n$, the set $\left\{\left(\xi(a_i)\right)_{i \leq k_n} :
  \sM \models P_n (\bar{a})\right\}$ is analytic.
\end{enumerate}
\end{dfn}

The function $\xi$ is called an $F$-parameterization of $\sM$.
Morozov also introduced the following stronger condition, essentially
requiring that $\sM$ be able to define its own
$F$-selfparameterization.

\begin{dfn}[\cite{fparam}]\label{fsparam} Let $\sM$ be an
  $F$-parameterizable structure.  We say that $\sM$ is weakly
  selfparameterizable if and only if there are functions $\Xi, p:
  \sM \times \omega \to \omega$, both definable without parameters in
  $HF(\sM)$, with the following properties:
\begin{enumerate}
\item For all $x \in \sM$ and all $m \in \omega$, we have $\Xi(x,m) =
  \xi(x)[m]$, and
\item For all $f \in \omega^\omega$ there is some $x \in \sM$ such
  that for all $n \in \omega$ we have $p(x,n) = f(n)$.
\end{enumerate}
\end{dfn}

In making sense of effectiveness on uncountable structures, a major
motivation is to describe a sense in which real number arithmetic ---
an operation that, while not Turing computable, does not seem horribly
ineffective --- can be considered to be effective.  Morozov proved the
following result, which is a major source of motivation for the
present paper.

\begin{prop}[\cite{fparam}]\label{realsarefparam} The real field is weakly
  $F$-selfparameterizable. \end{prop}

\begin{proof}[Outline of proof.] Define a function $\xi:\real \to
  \omega^\omega$ maps $x$ to its decimal expansion.  This function is
  definable without parameters in $HF(\real)$, in the sense
  required by Definition \ref{fsparam}.
\end{proof}

\subsection{$\Sigma$-definability}

In the present paper, when $\real$ denotes a structure, it is the
ordered field of real numbers.  The following definition is standard,
and appears in equivalent forms in \cite{barwisebk},
\cite{ershovopred}, and \cite{ershovhrm}.

\begin{dfn} Given a structure $\sM$ with universe $M$, we define a new
  structure $HF(\sM)$ as follows.
\begin{enumerate}
\item The universe of $HF(\sM)$ is the union of the chain
  $HF_n(M)$ defined as follows:
\begin{enumerate}
\item $HF_0(M) = \emptyset$
\item $HF_{n + 1}(M) = \sP^{< \omega}\left(M \cup
  HF_n(M)\right)$, where $\sP^{< \omega}(S)$ is the set of all finite
  subsets of $S$
\end{enumerate}
\item The language for $HF(\sM)$ consists of a binary predicate $\in$
  interpreted as membership, plus a predicate symbol $\sigma^*$ for
  each predicate symbol $\sigma$ of the language of $\sM$, interpreted
  in such a way that \[HF(\sM) \models \sigma^*(\{x_0\},\{x_1\}, \dots)
  \Leftrightarrow \sM \models \sigma(x_0, x_1, \dots),\] and similar
  additions for each function and constant symbol.
\end{enumerate}
\end{dfn}

Barwise \cite{barwisebk} noted the following connection to
computation, although he recorded that it was already well-known.

\begin{thm} Let $S$ be a relation on $\fN =
  \left(\mathbb{N}, +, \cdot\right)$.  Then 
\begin{enumerate}
\item $S$ is (classically) computably enumerable if and only if $S$ is
  $\Sigma_1$ on $HF(\fN)$.
\item $S$ is (classically) computable if and only if $S$ is $\Delta_1$
  on $HF(\fN)$.
\end{enumerate}
\end{thm}

This fact gives rise to Ershov's definition \cite{ershovopred} of a notion
generalizing computability to structures other than $\fN$.  We will
first give Barwise's definition of the class of $\Sigma$-formulas.

\begin{dfn} The class of $\Sigma$-formulas are defined by induction.
\begin{enumerate}
\item Each $\Delta_0$ formula is a $\Sigma$-formula.
\item If $\Phi$ and $\Psi$ are $\Sigma$-formulas, then so are
  $(\Phi \wedge \Psi)$ and $(\Phi \vee \Psi)$.
\item For each variable $x$ and each term $t$, if $\Phi$ is a
  $\Sigma$-formula, then the following are also $\Sigma$-formulas:
\begin{enumerate}
\item $\exists (x \in t) \ \  \Phi$
\item $\forall (x \in t) \ \ \Phi$, and
\item $\exists x \Phi$.
\end{enumerate}
\end{enumerate}
\end{dfn}

A predicate $S$ is called a \emph{$\Delta$-predicate} if both $S$ and
its complement are defined by $\Sigma$-formulas.

\begin{dfn} Let $\sM$ and $\sN = (N, P_0, P_1, \dots)$ be structures.
  We say that $\sN$ is $\Sigma$-definable in $HF(\sM)$ if and only if
  there are $\Sigma$-formulas $\Psi_0, \Psi_1, \Psi_1^*, \Phi_0,
  \Phi_0^*, \Phi_1, \Phi_1^*, \dots$ such that
\begin{enumerate}
\item $\Psi_0^{HF(\sM)} \subseteq HF(\sM)$ is nonempty,
\item $\Psi_1$ defines a congruence relation on $\left(\Psi_0^{HF(\sM)},
  \Phi_0^{HF(\sM)}, \Phi_1^{HF(\sM)}, \dots\right)$,
\item $(\Psi_1^{*})^{HF(\sM)}$ is the relative complement in
  $(\Psi_0^{HF(\sM)})^2$ of $\Psi_1^{HF(\sM)}$,
\item For each $i$, the set $(\Phi_i^{*})^{HF(\sM)}$ is the relative complement in
  $\Psi_0^{HF(\sM)}$ of $\Phi_i^{HF(\sM)}$, and
\item $\sN \simeq \left(\Psi_0^{HF(\sM)}, \Phi_0^{HF(\sM)},
  \Phi_1^{HF(\sM)}, \dots\right)/_{\Psi_1^{HF(\sM)}}$.
\end{enumerate}
\end{dfn}

\section{$\real$-computation and $F$-parameterizability}\label{compare}

We will first show that every $\real$-computable structure is
$F$-parameterizable.  We will also show that the converse of this
theorem does not hold.

\begin{thm} Let $\sM$ be a $\real$-computable structure.  Then $\sM$
  is $F$-parameterizable.\end{thm}

\begin{proof}
Let $\xi$ be an $F$-parameterization of the real field.  Then
consider the natural product map $\xi_\infty : \real_{\infty} \to
\omega^\omega$.  By Path Decomposition (Theorem 1 in Chapter 2 of
\cite{bcss}), the universe of $\sM$ must be a countable union of
semialgebraic sets of $\real_{\infty}$, and so $im (\xi_\infty \upharpoonright_{\sM})$ must
be analytic.  Similarly, for each predicate $P$ in the language of
$\sM$, the set $P(\sM)$ is a countable union of semialgebraic sets, so
that $\xi_{\infty}(P(\sM))$ is also analytic.  Thus, $\xi_\infty \upharpoonright_{\sM}$ is an
$F$-parameterization of $\sM$.
\end{proof}

The converse of this theorem is not true.  Indeed, a mathematically
familiar structure meets the stronger condition of being weakly
$F$-selfparameterizable but not $\real$-computable.

\begin{thm}\label{exp} The structure $(\real, +, \cdot, 0, 1, e^x)$ is weakly
  $F$-selfparameterizable but not $\real$-computable.\end{thm}

\begin{proof} Suppose $(\real, +, \cdot, 0, 1, f)$ is a
  $\real$-computable structure.  By Path Decomposition (Theorem 1 in
  Chapter 2 of \cite{bcss}), the atomic diagram must be a countable
  disjoint union of semi-algebraic sets, so that, in particular, the
  graph of $f$ must be a countable disjoint union of semi-algebraic
  sets.  The graph of $e^x$ does not have this property, so
  $(\real, +, \cdot, 0, 1, e^x)$ is not $\real$-computable.

Let $\xi$ be the standard weak $F$-selfparameterization of the reals.
It suffices to show that $\xi(e^x)$ is analytic.  We first
show that $\xi(\ln(x))$ is analytic, and the result for the
exponential function will follow.  Since \[\ln x = \int _1 ^x
\frac{1}{x} dx\] the actual value of $\ln x$ must be between the upper
and lower Riemann sum approximations to this integral.  If the
interval $[1,x]$ is partitioned into $n$ equal intervals, the upper
Riemann sum is \[ U_n = \sum\limits_{k=0}^{n-1} \left(\frac{x-1}{n}
\cdot \frac{1}{1+\frac{k(x-1)}{n}} \right) \] and the lower sum is \[
L_n = \sum\limits_{k=1}^{n} \left(\frac{x-1}{n} \cdot \frac{1}{1+\frac{k(x-1)}{n}} \right) \] so that we have the
error estimate \[ |\ln x - U_n| \leq U_n - L_n.\]  Now the difference
of the sums is given by \[U_n - L_n = \frac{x^2-2x+1}{nx}.\]
Consequently, we have the definition \[\ln x = y \Leftrightarrow
\bigwedge\limits_{n \in \omega} \hspace{-0.15in} \bigwedge |U_n - y|
\leq \frac{x^2 - 2x + 1}{nx},\] so that $\xi(\ln(x))$ is Borel.
This completes the proof.
\end{proof}

It seems natural to ask for some topological characterization of the
$\real$-computable structures among the $F$-parameterizable
structures.  We might hope, for instance, that by replacing
``analytic'' with ``$\mathbf{\Sigma^0_1}$'' or some such class in
Definition \ref{df:fparam}, we might find a class that coincides
exactly with the $\real$-computable structures.  A more careful
analysis of the foregoing proof and of Morozov's proof of Proposition
\ref{realsarefparam}, however, shows that such a characterization is
impossible, at least with the standard parameterization of the real field.

\begin{lem} Let $\xi$ be the weak $F$-selfparameterization of the real
  field given in the proof of Proposition \ref{realsarefparam}.  Then
  the following sets are $\mathbf{\Pi^0_1}$:
\begin{enumerate}
\item $\xi(=)$
\item $\xi(+)$
\item $\xi(\cdot)$
\end{enumerate}
\end{lem}

\begin{proof} For equality, it suffices to check, for each $n$,
  whether the decimal approximation of the two numbers up to $10^{-n}$
  either match exactly or match except for a terminal sequence of
  $9$'s starting at the $10^{-k}$ place, where the decimal
  approximation of the two numbers up to $10^{-k+1}$ differ by
  $10^{-k+1}$.  Since this is a $\mathbf{\Delta^0_1}$ condition for
  each $n \in \omega$, the equality relation is $\mathbf{\Pi^0_1}$.

The cases of addition and multiplication are nearly identical to one
another.  Let $x_n$ be the decimal approximation of $x$ up to
$10^{-n}$.  To determine whether $x + y = z$, it suffices to check
whether, for each $n$, we have \[|x_n + y_n - z_n| \leq 10^{-n+1}.\]
Again, the condition is $\mathbf{\Delta^0_1}$ for each $n \in \omega$, so that
the addition relation is $\mathbf{\Pi^0_1}$.
\end{proof}

\begin{lem} The set $\xi(=)$ is $\mathbf{\Pi^0_1}$ hard.\end{lem}

\begin{proof}
Let $Q$ be a $\mathbf{\Pi^0_1}$ set, defined by $\forall x
\tilde{Q}(x,S)$, where $S$ ranges over sets of natural numbers.  For
each set $S$, define the set $T_S$ to be the set of all $x$ such that
$\lnot \tilde{Q}(x,S)$.  Now $Q(S)$ holds if and only if $T(S)$ is
empty.

Given a set $S \subseteq \omega$, we can view $S$ as a real number in
$[0,1]$ in a natural way, as \[x_S = \sum\limits_{i = 0}^\infty
\chi_S(i) 10^{-i}.\]  Now to decide $Q(S)$ from $\xi(=)$, it suffices
to check whether $\xi(=)$ holds of the pair $(\xi(x_{T_S}),
\xi(x_\emptyset))$.
\end{proof}

Now in the proof of Theorem \ref{exp}, we saw that $\xi(e^x)$ is $\mathbf{\Pi^0_1}$ over
$\xi(=) \oplus \xi(+) \oplus \xi(\cdot)$, but this set is
$\mathbf{\Pi^0_1}$ complete.  Consequently, in the standard
parameterization of the real field, the Borel hierarchy does
not distinguish between the complexity of the field structure and that
of the exponential function.  Of course, it may be
that there is a simpler parameterization of the real field in which
this proof does not suffice, so we cannot as yet categorically rule
out any such characterization.

\section{$\real$-computation and $\Sigma$-definability}\label{rcmtsig}

It should be noted that the work of this section is related to earlier
work of Ashaev, Belyaev, and Myasnikov \cite{ashaevbelyaevmyasnikov},
who proved a similar but weaker result for the much more general
definition of computation over an arbitrary structure given by Blum,
Shub, and Smale.

\subsection{$\real$-Computable Structures are $\Sigma$-definable in $HF(\real)$}\label{rcsd}

\begin{prop} Every $\real$-computable structure is $\Sigma$-definable
  in $HF(\real)$.\end{prop}

\begin{proof} Let $\sM = (M, P_0, P_1, \dots)$ be a $\real$-computable
  structure.  Now $M$, and each $P_i^\sM$, are $\real$-computable
  sets, and so, in particular, both they and their complements are
  $\real$-semidecidable.

\begin{lem} The class $\mathcal{RC}$ of $\real$-computable functions
  is the smallest class containing all polynomial functions, the
  characteristic function of $<$, and the shift function, and closed
  under composition, juxtaposition, primitive recursion, and
  minimalization. \end{lem}

\begin{proof} The finite-dimensional case of this result, including
  most of the important issues, is given in \cite{bss}.  The
  polynomial functions, the characteristic function of $<$, and the
  shift function are clearly $\real$-computable.  That the class of
  $\real$-computable functions is closed under composition,
  juxtaposition, primitive recursion, and minimalization is shown in
  \cite{bss}.  On the other hand, let $f$ be a $\real$-computable
  function computed by a machine $\Omega_f$.  Then $\Omega_f$
  describes how to build $f$ using only composition, juxtaposition,
  primitive recursion, and minimalization, from polynomials, the
  characteristic function of $<$, and the shift function.\end{proof}

Now in $HF(\real)$, we will represent elements $(x_i)_{i \in \mathbb{Z}} \in
\real_\infty$ as sets of the form
$\left\{\left\{\{i\},\{\{x_i\}\}\right\} : i \in I\right\}$, where
$I \subseteq \mathbb{Z}$ is the (finite) set of indices for nonzero elements.

\begin{lem}\label{rcompsdeffn} Every $\real$-computable function $f$ is a $\Sigma$-function
  in $HF(\real)$.  Moreover, the complement of the graph of $f$ is
  also defined by a $\Sigma$ formula.\end{lem}

\begin{proof} The polynomial functions and the characteristic function of $<$, by
    the definition of the structure $HF(\real)$, are definable without
    quantifiers in $HF(\real)$, as are their complements.  The shift
    operators are defined as
    \[\sigma_\ell\left(\left\{\left\{\{i\},\{\{x_i\}\}\right\} : i \in
    I\right\}\right) = \left(\left\{\left\{\{i\},\{\{y_i\}\}\right\} : i \in
    I\right\}\right) \Leftrightarrow \bigwedge\limits_{i \in I} y_i =
    x_{i+1}\]
and 
\[\sigma_r\left(\left\{\left\{\{i\},\{\{x_i\}\}\right\} : i \in
    I\right\}\right) = \left(\left\{\left\{\{i\},\{\{y_i\}\}\right\} : i \in
    I\right\}\right) \Leftrightarrow \bigwedge\limits_{i \in I} y_i =
    x_{i-1},\]
both of which can be expressed by quantifier-free formulas in
    $HF(\real)$, as can their complements.  The juxtaposition of two
    $\Sigma$ functions is a $\Sigma$-function  (by clause 2 of the
    definition of $\Sigma$-formulas), and if the complements of their
    graphs are $\Sigma$-definable, then the same is true of the
    juxtaposition.  The case of the complement is the same (by
    clause 3c).

Suppose that $F: \mathbb{Z}^{\geq 0} \times \real^\ell \to \real$ is a
$\Sigma$-function (defined, say, by $\Phi(t,\bar{x}, y)$) and that the
complement of its graph is $\Sigma$-definable (by $\overline{\Phi}(t,
\bar{x}, y)$).  Then $L(\bar{x}) := \mu t[F(t,\bar{x}) = 0]$ is
defined by \[\Theta(\bar{x}, t) := \Phi(t, \bar{x}, 0) \wedge \forall
s\in t [\lnot \Phi(s, \bar{x}, 0)],\] and its complement is defined
similarly.

To show that the class of $\Sigma$-definable functions is closed under
primitive recursion, suppose that $h$ and $g$ are $\Sigma$-definable
functions, and $f$ is defined by the following schema:
\begin{eqnarray*} f(0,x_2,\dots, x_k)&=&g(x_2, \dots, x_k)\\f(y+1, x_2,
  \dots, x_k)&=&h(y,f(y, x_2, \dots, x_k).\end{eqnarray*}
Now since $g$ and $h$ are $\Delta$ relations, we can apply $\Delta$
recursion (Corollary I.6.6 of \cite{barwisebk}) to show that $f$ is
also a $\Delta$ relation.
\end{proof}

Now since the characteristic function of $M$, and of each $P_i$ is
$\real$-computable, they are also $\Sigma$ functions of $HF(\real)$,
as are the characteristic functions of their complements.  Let us say
that, as $U$ ranges over $\{M, M^c, P_0, P_0^c, P_1, P_1^c, \dots)$,
the function $\chi_U$ is defined by $\Theta_U(x, \{y\})$; that is,
suppose that $\chi_U(x) = y$ if and only if $\Theta_U(x,\{y\})$.  Now
$\Theta_U(x,\{1\})$ is a $\Sigma$-formula defining $U$.  We let
$\Psi_1(x,y)$ be the relation $x=y$, so that both $\Psi_1$ and its
complement are defined by $\Sigma$ formulas.  This completes the
proof.
\end{proof}

R.\ Miller and Mulcahey \cite{rmiller2} raise the issue, in the
context of a different notion of effectiveness for uncountable
structures, of ``simulating'' an uncountable structure by a
classically computable one.  The following result shows a way in which
this goal can be realized for certain $\real$-computable structures.

\begin{dfn} Let $\sA$ and $\sB$ be structures in a common signature.
  We write that $\sA \leq_1 \sB$ if $\sA$ is a substructure of $\sB$,
  and for all existential formulas $\varphi(\bar{x})$ and for all
  tuples $\bar{a} \subseteq \sA$, we have \[\sB \models
  (\varphi(\bar{a}) \Rightarrow \sA \models
  \varphi(\bar{a}).\]\end{dfn}

\begin{cor}\label{capprox} For any $\real$-computable structure $\sM$
  whose defining machine involves only computable reals as parameters,
  there is a computable structure $\sM^*$ such that $\sM^* \leq_1
  \sM$.\end{cor}

\begin{proof} Morozov and Korovina \cite{morozovkorovina} proved this
  result under the alternate hypothesis that $\sM$ is
  $\Sigma$-definable without parameters over $HF(\real)$.  In fact, if
  the field of computable real numbers is used in place of the field
  of real algebraic numbers, their proof works without further
  modification for a structure $\Sigma$-definable with finitely many
  computable parameters.  Since, by the theorem, $\sM$ satisfies this
  hypothesis, the result follows.\end{proof}

\subsection{Structures $\Sigma$-Definable in $HF(\real)$ are
  $\real$-Computable}\label{sdrc}

\begin{thm} Every structure $\Sigma$-definable over $HF(\real)$ has an
  isomorphic copy which is $\real$-computable.\end{thm}

\begin{proof} Let $\sM = (M, P_0, P_1, \dots)$ be $\Sigma$-definable
  over $HF(\real)$, via the scheme \[\Psi_0, \Psi_1, \Psi_1^*, \Phi_0,
  \Phi_0^*, \Phi_1, \Phi_1^*, \dots.\]  We may assume that all of these
  formulae are in prenex normal form.

\begin{lem}\label{enumhf} There is an enumeration function $e:\real_\infty \to HF(\real)$ with the
  following properties:
\begin{enumerate}
\item $e$ is a bijection;
\item $e$ is a $\Sigma$ function on $HF(\real)$;
\item $e$ is $\real$-computable; and
\item For any $\Delta_0$ formula $\varphi(\bar{x})$, the relation
  $HF(\real) \models \varphi(e(n_1), \dots, e(n_k))$ is
  $\real$-computable.
\end{enumerate}
\end{lem}

\begin{proof} 

Let $\sT$ be the class of finite trees, and let $T: \omega \to \omega$
be a classically computable Friedberg enumeration of $\sT$.  That such
an enumeration must exist is guaranteed by results in \cite{cckm}.
Let $T_k : \omega \to \omega$ be a classically computable Friedberg
enumeration of the elements of $\sT$ with exactly $k$ terminal nodes.
Let $\epsilon: \mathbb{Z} \to \omega$ by \[\epsilon (i) :=
\left\{\begin{array}{ll} 2i & \mbox{if $i \geq 0$}\\ -2i-1 & \mbox{if
  $i < 0$}\end{array}\right. .\]

Let $\mathbf{x} := (x_i : i \in \mathbb{Z}) \in \real_\infty$.
Suppose $i_0$ is the first and $i_0+k$ the last such that $x_{i_0}$
and $x_{i_0+k}$ are nonzero.  Now to specify $e(\mathbf{x})$, we take
the classically computable tree with index $T_{k+1}(\epsilon(i_0))$,
order its terminal nodes lexicographically from $0$ to $k$, and label
terminal node $j$ with the real number $x_j$.  The interpretation of
this tree with these labels as a member of $HF(\real)$ is
straightforward.  

The function $e$ defined in this way is bijective and
$\real$-computable.  By Lemma \ref{rcompsdeffn}, it is also a
$\Sigma$-function in $HF(\real)$.  To decide if $e(n) \in e(m)$, it
suffices to check whether there is an embedding of the appropriate
finite trees, preserving the real labels on terminal nodes.  This
operation is $\real$-computable.  Satisfaction of other $\Delta_0$
formulas is $\real$-computable, by induction on their form.
\end{proof}

Let $U$ range over the various $\Psi_i, \Psi_i^*, \Phi_i$, and
$\Phi_i^*$.  Now suppose that $U$ is of the form $\exists y_1, \dots,
y_n \varphi(\bar{y}, \bar{x})$, where $\varphi$ is $\Delta_0$.  Then
$U(\bar{x})$ holds exactly when \[\exists j_1, \dots, j_n, \ell_1,
\dots, \ell_m \left[\left(\bigwedge\limits_{i \leq m} e(\ell_i) =
  x_i\right) \wedge \varphi(e(j_1), \dots, e(j_n), e(\ell_1), \dots,
  e(\ell_m)\right].\]  The part within the brackets is
$\real$-computable.  

Now suppose both $U$ and its complement have definitions of the form
$\exists \bar{x} [R(\bar{x})]$, where $R$ is $\real$-computable.  Using Path
Decomposition \cite{bcss}, we can write $R(\bar{x})$ as a countable disjoint
union of semialgebraic sets $R_i(\bar{x})$, so that $U$ is defined by
\[\exists \bar{x} \bigwedge\limits_{i \in
  \omega}\hspace{-0.15in}\bigwedge R_i(\bar{x}) \Leftrightarrow \bigwedge\limits_{i \in
  \omega}\hspace{-0.15in}\bigwedge \exists \bar{x} R_i(\bar{x}).\]  By
  Tarski's elimination theory for the real ordered field, the
  condition on the right is a countable disjoint union of
  semialgebraic sets.  Now both $\Psi_1$ and each of the $\Phi_i$ are
  $\real$-computable. 

It remains to show that $\sM$ has an isomorphic copy in which $\Psi_0$
is $\real$-computable.  We may suppose $\Psi_0$ is in prenex normal
form, so that $\Psi_0(\bar{x}) = \exists \bar{y} [\psi(\bar{x},
  \bar{y})]$, where $\psi$ is quantifier-free.  Let $M' \subseteq
(HF(\real)^{n_0} \times \real_\infty)$, where $n_0$ is the arity of
$\Psi_0$, be such that $(\bar{x}, \mathbf{t}) \in M'$ if and only if
$\psi(\bar{x}, e(\mathbf{t}))$.  We will interpret all symbols of the
language by their meaning on the first coordinate (so that, in
particular, $(M', \Phi_0, \Phi_1, \dots) \models (\bar{x}, \mathbf{t})
= (\bar{x}, \mathbf{s})$ for all $\mathbf{t}, \mathbf{s}$.  Since $\psi$ is
quantifier-free, the structure $(M', \Phi_0, \Phi_1, \dots)$ is
$\real$-computable.
\end{proof}

\bibliographystyle{amsplain}
\bibliography{rcmt}

\end{document}